\documentclass[12pt,reqno]{amsart}
\usepackage{hyperref}
\usepackage{cleveref}
\usepackage{amssymb, eqnarray, mathrsfs, xcolor}
\usepackage{tikz}
\usetikzlibrary{arrows.meta, positioning}
 
\theoremstyle{plain}
\newtheorem{theorem}{Theorem}[section]

\newtheorem{lemma}[theorem]{Lemma}

\theoremstyle{definition}

\newcommand{\Pelc}[1]{\left(\Sigma X_n\right)_p}

\DeclareSymbolFont{bbold}{U}{bbold}{m}{n}
\DeclareSymbolFontAlphabet{\mathbbold}{bbold}

\DeclareMathOperator{\rank}{rank}

\newcommand{\cN}{\mathcal{N}}

\begin{document}

\title[Positive commutators]{Positive commutators of positive square-zero operators} 

\author{Roman Drnov\v sek}
\address{Faculty of Mathematics and Physics, University of Ljubljana,
  Jadranska 19, 1000 Ljubljana, Slovenia \ \ and \ \ \
  Institute of Mathematics, Physics, and Mechanics, Jadranska 19, 1000 Ljubljana, Slovenia}
\email{roman.drnovsek@fmf.uni-lj.si}

\author{Marko Kandi\'c}
\address{Faculty of Mathematics and Physics, University of Ljubljana,
  Jadranska 19, 1000 Ljubljana, Slovenia \ \ and  \ \ \
  Institute of Mathematics, Physics, and Mechanics, Jadranska 19, 1000 Ljubljana, Slovenia}
\email{marko.kandic@fmf.uni-lj.si}


\date{\today}
\baselineskip 7mm

\begin{abstract}
In this paper we first consider the question which nonnegative matrices are commutators of nonnegative square-zero matrices. Then, we treat infinite-dimensional analogues of these results for operators on the Banach lattices $L^p[0,1]$ and $\ell^p$ ($ 1 \leq p < \infty$). In the last setting we need to extend the notion of the nonnegative rank of a nonnegative matrix.
\end{abstract}

\maketitle

\section{Introduction}

Positive commutators of positive operators on Banach lattices have been the subject of extensive research; see e.g. \cite{Bracic:10, DK11, Drnovsek:11, Gao:14, KS17, KS17b, Drnovsek:18,DK19, DK25a, DK25b}. A systematic investigation of their properties began with \cite{Bracic:10}, where the authors studied the spectral properties of
the positive commutator $[A, B] := AB - BA$ formed by positive compact operators $A$ and $B$.
In particular, they showed that if $A$ and $B$ are nonnegative matrices such that the commutator $C=[A, B]$
is nonnegative as well, then, up to similarity with a permutation matrix, $C$ is a strictly upper triangular matrix, and so it is nilpotent. Recently, the authors of \cite{MRZ25} have studied commutators of square-zero operators on Hilbert spaces. In particular, they completely characterized commutators of square-zero matrices
by proving that commutators are, up to similarity, precisely direct sums $A \oplus (-A) \oplus N^2$, where 
$A$ is an invertible matrix and $N$ is a nilpotent matrix.

In this paper we first study the question which nonnegative nilpotent matrices are commutators of nonnegative square-zero matrices. This is contained in Section 2. Infinite-dimensional analogues of these results for operators on the Banach lattices $L^p[0,1]$ and $\ell^p$ ($ 1 \leq p < \infty$) are considered in Section 3, where  the nonnegative rank of positive operators on vector lattices is also introduced.

Since some of our results hold in general setting of vector lattices, we recall some basic
definitions and properties of vector lattices and operators on them. For the terminology
and details not explained here we refer the reader to \cite{LZ71} or \cite{AB06} or \cite{AA02}.
Let $X$ be a vector lattice with the positive cone $X^+$. 
A subspace $A$ of $X$ is called an {\it ideal} whenever $|x| \leq |y|$ and $y \in A$ imply $x \in A$.
In the finite-dimensional case $X = \mathbb R^n$ ideals are also called {\it standard subspaces}, since they are precisely linear spans of some standard unit vectors.
An order dense ideal of $X$ is said to be a {\it band}. In the case $X=L^p(\mu)$ ($1 \leq p < \infty$) bands are precisely (norm) closed ideals.
The band
$$ S^d := \{ x \in X : |x| \land |y| = 0 \ \ {\rm for \ all} \  y \in S \}  $$
is called the {\it disjoint complement} of a set $S$ of $X$.

A linear operator $T$ on $X$ is {\it positive} if $T$ leaves the positive cone $X^+$ invariant, that is, $T(X^+) \subseteq X^+$.
In the finite-dimensional case $X = \mathbb R^n$ these are precisely nonnegative $n \times n$ matrices. 
Let $T$ be a positive operator on a vector lattice $X$. The {\it null ideal} $\cN(T)$ is the ideal in $X$ defined by
$$ \cN(T) = \{ x \in X : T |x| = 0 \} .  $$

The {\it nonnegative rank} of a nonnegative $m \times n$ matrix $A$ is equal to the smallest nonnegative integer $k$ for which there exist a nonnegative $m \times k$ matrix $L$ and a nonnegative $k \times n$ matrix $R$ such that $A = L R$. 
We denote it by $\rank^+ (A)$. 
To obtain the usual rank (denoted by $\rank (A))$, we drop the condition that $L$ and $R$ must be nonnegative. 
It is easy to see that $\rank^+ (A)$ is the smallest nonnegative integer $k$ such that there exist nonnegative vectors $u_1$, $\ldots$, $u_k$ and 
$v_1$, $\ldots$, $v_k$ such that $A = \sum_{i=1}^k u_i v_i^T$.

Finally, we recall Bohnenblust's result. By \cite[Theorem 7.1]{Bohnenblust}, every separable infinite-dimensional Banach lattice $L^p(\mu)$ $(1\leq p<\infty)$ is isometric and order isomorphic to one of the following Banach lattices: $\ell^p$, $L^p[0,1]$, $\ell^p\oplus L^p[0,1]$ or $\ell_n^p\oplus L^p[0,1]$.

\section{The finite-dimensional case}

In this section we study the question which nonnegative matrices are commutators of nonnegative square-zero matrices. We begin with the special case.

\begin{theorem}
\label{main-finite_dimensional}
\label{equivalent}
     For a nonnegative $n \times n$ matrix $T$ the following assertions are equivalent. 
     \begin{enumerate}
        \item There exist nonnegative $n \times n$ matrices $M$ and $N$ such that $T=M N$ and 
        $M^2=N^2=N M=0$. 
        \item There exists a nonnegative $n \times n$ matrix $U$ such that $T=U^2$ and $U^3=0$.
        \item There exists a decomposition $\mathbb R^n = L_1\oplus L_2\oplus L_3$ on standard subspaces with respect to which the operator $T$ is of the form 
        $$T=\left(\begin{matrix}
             0 & 0 & T_{13}\\
             0 & 0 & 0 \\
             0 & 0 & 0
        \end{matrix}\right)$$ 
        for some positive operator $T_{13}\colon L_3\to L_1$ with  $\rank^+ (T_{13}) \leq \dim (L_2)$.
     \end{enumerate}
\end{theorem}

\begin{proof}
(i)$\Rightarrow$(ii)  Define $U:= M+N$. Since $M^2=N^2=N M=0$, we have $U^2 = M N = T$ and 
$U^3 = M N (M+N) = 0$. 

(ii)$\Rightarrow$(iii) Define standard subspaces $L_1=\mathcal N(U)$, $L_2=\mathcal N(U)^d\cap \mathcal N(U^2)$ and $L_3=\mathcal N(U^2)^d$. 
 With respect to the decomposition $\mathbb R^n = L_1\oplus L_2\oplus L_3$ the positive operators $U$ and $U^2$ are of the form 
    $$U=\left(\begin{matrix}
    0 & U_{12} & U_{13}\\
    0 & 0 & U_{23} \\
    0 & 0 & 0
        \end{matrix}\right) \qquad \textrm{and} \qquad T=U^2=\left(\begin{matrix}
    0 &  0 & U_{12} U_{23}\\
    0 & 0 & 0 \\
    0 & 0 & 0
        \end{matrix}\right).$$
Since $T_{13} = U_{12} U_{23}$, we have $\rank^+ (T_{13}) \leq \dim (L_2)$.

(iii)$\Rightarrow$(i) Denote $d_i = \dim (L_i)$ for $i=1, 2, 3$. Since $\rank^+ (T_{13}) \leq d_2$,
there exist a nonnegative $d_1 \times d_2$ matrix $L$ and a nonnegative $d_2 \times d_3$ matrix $R$ such that $T_{13} = L R$.  
Define the nonnegative $n \times n$ matrices by
$$ M = \left(\begin{matrix}
    0 & L & 0 \\
    0 & 0 & 0 \\
    0 & 0 & 0
        \end{matrix}\right) \qquad \textrm{and} \qquad N=\left(\begin{matrix}
    0 & 0 & 0\\
    0 & 0 & R \\
    0 & 0 & 0
        \end{matrix}\right).$$
Then it is easy to see that $T=M N$ and $M^2=N^2=N M=0$. 
\end{proof}

We now find necessary conditions for a nonnegative matrix to be a commutator of nonnegative square-zero matrices.

\begin{theorem}
\label{necessary}
    Let $T$ be a nonnegative $n \times n$ matrix such that $T= M N - N M$ for some 
    nonnegative $n \times n$ matrices $M$ and $N$ with $M^2=N^2=0$. 
    Then $M T = T M = N T = T N = 0$ and 
    there exists a decomposition $\mathbb R^n = L_1\oplus L_2\oplus L_3$ on standard subspaces with respect to which the operator $T$ is of the form 
        $$T=\left(\begin{matrix}
             0 & 0 & T_{13}\\
             0 & 0 & 0 \\
             0 & 0 & 0
        \end{matrix}\right)$$ 
        for some positive operator $T_{13}\colon L_3\to L_1$ with $\rank (T_{13}) \leq \dim (L_2)$.
\end{theorem}

\begin{proof}
Since $M^2 = 0$, we have $0 \leq M T = - M N M \leq 0$, and so $M T = 0 = M N M$.
Similarly, since $N^2 = 0$, we have $0 \leq T N = - N M N \leq 0$, and so $T N = 0 = N M N$.
It follows that $T M = M N M = 0$ and $N T = N M N = 0$. 

Denote $S = M+N$. Since $S^2 = M N + N M$, we have $S^3 = (M N + N M)(M + N) = 0$.
Now, define standard subspaces $L_1=\mathcal N(S)$, $L_2=\mathcal N(S)^d\cap \mathcal N(S^2)$ and $L_3=\mathcal N(S^2)^d$. 
With respect to the decomposition $\mathbb R^n = L_1\oplus L_2\oplus L_3$ the positive operators $S$ and $S^2$ are of the form 
    $$S=\left(\begin{matrix}
    0 & \star & \star \\
    0 & 0 & \star \\
    0 & 0 & 0
        \end{matrix}\right) \qquad \textrm{and} \qquad S^2=\left(\begin{matrix}
    0 & 0 & \star \\
    0 & 0 & 0 \\
    0 & 0 & 0
        \end{matrix}\right).$$
Since $0 \leq T = M N - N M \leq M N + N M = S^2$, the operator $T$ has the same pattern as the operator $S^2$, so that 
    $$T=\left(\begin{matrix}
             0 & 0 & T_{13}\\
             0 & 0 & 0 \\
             0 & 0 & 0
        \end{matrix}\right)$$ 
        for some positive operator $T_{13}\colon L_3\to L_1$. 
Since $0 \leq M, N \leq S$, the operators $M$ and $N$ have the following forms
$$M=\left(\begin{matrix}
    0 & M_{12} & M_{13}\\
    0 & 0 & M_{23} \\
    0 & 0 & 0
        \end{matrix}\right) \qquad \textrm{and} \qquad N =\left(\begin{matrix}
    0 & N_{12} & N_{13}\\
    0 & 0 & N_{23} \\
    0 & 0 & 0
        \end{matrix}\right).$$
Since 
$$ T = M N - N M=\left(\begin{matrix}
    0 & 0 & M_{12} N_{23} - N_{12} M_{23}\\
    0 & 0 & 0 \\
    0 & 0 & 0
        \end{matrix}\right) , $$
it follows that $T_{1 3} = M_{12} N_{23} - N_{12} M_{23}$. Since $M^2 = N^2 = 0$, we have 
$M_{12} M_{23} =0$ and $N_{12} N_{23} = 0$, and so we can write 
$T_{13} = (M_{12} - N_{12})(M_{23}+N_{23})$. This implies that $\rank (T_{13}) \leq \dim (L_2)$ as desired.
\end{proof}

In view of \Cref{equivalent}, we would expect that the converse implication in \Cref{necessary} also holds. However, we will show in \Cref{example} that this is not the case.

Define the $4 \times 4$ matrix
\begin{equation}
T_{1 3} =\left(\begin{matrix}
             1 & 1 & 0 & 0 \\
             1 & 0 & 1 & 0 \\
             0 & 1 & 0 & 1 \\
             0 & 0 & 1 & 1 
        \end{matrix}\right) . 
\label{T13}        
\end{equation}

Note that $\rank (T_{1 3})=3$ and $\rank^+ (T_{13})=4$; see e.g. \cite{CR93}. In the proof of \Cref{example} we will need the following lemma.

\begin{lemma}\label{rank3}
Let $u =(u_1 \ u_2 \ u_3 \ u_4)^T$ and $v =(v_1 \ v_2 \ v_3 \ v_4)^T$ be nonnegative (column) vectors. Then 
$$ \rank(T_{1 3} + u v^T) \geq 3 . $$
\end{lemma}

\begin{proof}
    We distinguish 2 cases. 
    
    {\it Case 1}: $u_1 \leq u_3$. We consider the submatrix obtained from $T_{1 3} + u v^T$ by deleting the first column and the last row, that is the matrix
    $$ S := \left(\begin{matrix}
             1 & 0 & 0 \\
             0 & 1 & 0 \\
             1 & 0 & 1 
        \end{matrix}\right) + 
        \left(\begin{matrix}
             u_1 \\
             u_2 \\ 
             u_3 
           \end{matrix}\right) 
         \left(\begin{matrix}
             v_2 & v_3 & v_4
           \end{matrix}\right) . $$
Denoting $$ A = \left(\begin{matrix}
             1 & 0 & 0 \\
             0 & 1 & 0 \\
             1 & 0 & 1 
        \end{matrix}\right) , \ \tilde{u} = \left(\begin{matrix}
             u_1 \\
             u_2 \\ 
             u_3 
           \end{matrix}\right) \qquad \textrm{and} \qquad \tilde{v} = \left(\begin{matrix}
             v_2 \\
             v_3 \\
             v_4
           \end{matrix}\right) ,  $$
by the matrix determinant lemma (see e.g. \cite{Wiki} or \cite[Lemma 1.1]{DZ07}), the determinant of $S$ is equal to 
$$ (1 + \tilde{v}^T A^{-1} \tilde{u}) \, {\rm det}(A) = 1 + \left(\begin{matrix}
             v_2 & v_3 & v_4
           \end{matrix}\right)
           \left(\begin{matrix}
             1 & 0 & 0 \\
             0 & 1 & 0 \\
             -1 & 0 & 1 
        \end{matrix}\right)           
           \left(\begin{matrix}
             u_1 \\
             u_2 \\ 
             u_3 
           \end{matrix}\right) = $$ 
$$ = 1 + v_2 u_1 + v_3 u_2 + v_4 (u_3-u_1) \geq 1 . $$
This proves that the rank of $T_{1 3} + u v^T$ is at least $3$.

{\it Case 2}: $u_3 \leq u_1$. We now consider the submatrix obtained from $T_{1 3} + u v^T$ by deleting the second row and the last  column, that is the matrix
    $$ R := \left(\begin{matrix}
             1 & 1 & 0 \\
             0 & 1 & 0 \\
             0 & 0 & 1 
        \end{matrix}\right) + 
        \left(\begin{matrix}
             u_1 \\
             u_3 \\ 
             u_4 
           \end{matrix}\right) 
         \left(\begin{matrix}
             v_1 & v_2 & v_3
           \end{matrix}\right)  . $$
By the matrix determinant lemma again, the determinant of $R$ is equal to 
$$ 1 + \left(\begin{matrix}
             v_1 & v_2 & v_3
           \end{matrix}\right)
           \left(\begin{matrix}
             1 & -1 & 0 \\
             0 & 1 & 0 \\
             0 & 0 & 1 
        \end{matrix}\right)           
           \left(\begin{matrix}
             u_1 \\
             u_3 \\ 
             u_4 
           \end{matrix}\right) = 1 + v_1 (u_1-u_3) + v_2 u_3 + v_3 u_4 \geq 1 , $$
showing that the rank of $T_{1 3} + u v^T$ is at least $3$.
\end{proof}

We now show that the necessary conditions of \Cref{necessary} are not sufficient.

\begin{theorem}
Define the matrix 
$$T=\left(\begin{matrix}
             0 & 0 & T_{13}\\
             0 & 0 & 0 \\
             0 & 0 & 0
        \end{matrix}\right)  $$ 
with respect to the decomposition $\mathbb R^{11} = \mathbb R^4 \oplus \mathbb R^3 \oplus \mathbb R^4$,
where the matrix $T_{13}$ is given by (\ref{T13}). Then $T$ cannot be written as a commutator $M N - N M$, where $M$ and $N$ are 
nonnegative $11 \times 11$ matrices with $M^2 = N^2 =0$.
\label{example}
\end{theorem}

\begin{proof}
Observe that $\rank (T_{1 3})=3$, and so for $T$ the necessary condition of \Cref{necessary} is satisfied,  while $\rank^+ (T_{13})=4$, so that the assertion (iii) of \Cref{equivalent} is not true.
Assume to the contrary that $T = M N - N M$, where $M$ and $N$ are nonnegative $11 \times 11$ matrices with $M^2 = N^2 =0$.
Since $M T = T M = N T = T N = 0$ by \Cref{necessary} and $T_{1 3}$ has no zero rows and columns, $M$ and $N$ must have the forms
$$M=\left(\begin{matrix}
    0 & M_{12} & M_{13} \\
    0 & M_{22} & M_{23} \\
    0 & 0 & 0
        \end{matrix}\right) \qquad \textrm{and} \qquad N =\left(\begin{matrix}
    0 & N_{12} & N_{13} \\
    0 & N_{22} & N_{23} \\
    0 & 0 & 0
    \end{matrix}\right) $$
with respect to the decomposition $\mathbb R^{11} = \mathbb R^4 \oplus \mathbb R^3 \oplus \mathbb R^4$.
Since $T = M N - N M$, we have $T_{13} = M_{12} N_{23} - N_{12} M_{23}$ or 
$T_{13} + N_{12} M_{23} = M_{12} N_{23}$. 
Since $T_{13} \neq 0$, it follows that $M_{12} N_{23} \neq 0$, 
so that $M_{12} \neq 0$. Since $\rank^+ (T_{13})=4$, it must hold that 
$N_{12} M_{23} \neq 0$, so that $M_{23} \neq 0$. 
From $M^2 = 0$ it follows that $M_{12} M_{23} = 0$, and so 
$1 \leq \rank (M_{12}) \leq 2$ and $1 \leq \rank (M_{23}) \leq 2$.
We distinguish 2 cases.

{\it Case 1}: $\rank(N_{12} M_{23}) = 1$. Then $N_{12} M_{23} = u v^T$ 
for some nonnegative vectors $u$ and $v$. It follows that 
$\rank(T_{1 3} + u v^T) = \rank(M_{12} N_{23}) \leq \rank(M_{12}) \leq 2$. This is a contradiction with the conclusion of Lemma \ref{rank3}.

{\it Case 2}: $\rank(N_{12} M_{23}) = 2$. Then $\rank^+ (N_{12} M_{23}) = 2$ 
by \cite[Theorem 4.1]{CR93}, and so there exist nonnegative vectors $u$, $v$, $w$ and $z$ such that $N_{12} M_{23} = u v^T + w z^T$.
Now, $\rank(M_{23}) = 2$, and so $M_{23}$ has at least $2$ non-zero rows. 
Since $M_{12} M_{23} = 0$, $M_{12}$ is a $4 \times 3$ matrix and $M_{23}$ is a $3 \times 4$ matrix, 
we conclude that $M_{12}$ has only one non-zero column, so that $\rank(M_{12}) = 1$.
It follows that 
$$ \rank(T_{1 3} + u v^T) = \rank(M_{12} N_{23} - w z^T) \leq \rank(M_{12} N_{23}) + 1 
\leq \rank(M_{12}) + 1 = 2 . $$
This is again a contradiction with Lemma \ref{rank3}.
\end{proof}

\section{The infinite-dimensional case}

In this section we consider infinite-dimensional analogues of \Cref{main-finite_dimensional}. 
We first treat positive operators on the atomic Banach lattice $\ell^p$ $(1\leq p<\infty)$. 
In this setting we need to extend the notion of a nonnegative rank of a nonnegative matrix.

A positive operator $T\colon X\to Y$ between Archimedean vector lattices is said to have a \emph{finite nonnegative rank} if there exists a finite-dimensional Archimedean vector lattice $Z$ and positive operators $L\colon Z\to Y$ and $R\colon X\to Z$ such that $T= LR$. 

\[
\begin{tikzpicture}[>=Stealth, node distance=2.5cm, on grid]
  \node (X) at (0,0) {$X$};
  \node (Y) at (4,0) {$Y$};
  \node (Z) at (2,-2) {$Z$};

  \draw[->] (X) -- (Y) node[midway, above] {$T$};
  \draw[->] (X) -- (Z) node[midway, left, xshift=-4pt, yshift=-2pt] {$R$};
  \draw[->] (Z) -- (Y) node[midway, right, xshift=4pt, yshift=-2pt] {$L$};
\end{tikzpicture}
\]

Since every finite-dimensional Archimedean vector lattice is lattice isomorphic to $\mathbb R^k$ (ordered coordinatewise) for some $k\in \mathbb N_0$, in the definition of a finite nonnegative rank we can replace $Z$ with $\mathbb R^{\dim Z}$. 
Therefore, $T$ has a finite nonnegative rank whenever it factors through $\mathbb R^k$ (ordered coordinatewise) with positive factors. The \emph{nonnegative rank} $\rank^+(T)$ of $T$ is defined as the minimal $k\in \mathbb N_0$ above such that $T$ factors through $\mathbb R^k$ with positive factors. Therefore, this definition, in the case when $X=\mathbb R^n$ and $Y=\mathbb R^m$ coincides with the definition of a nonnegative rank of nonnegative $m\times n$ matrices.  
We define $\rank^+(T)=\infty$ if $T$ does not admit a factorization through any finite-dimensional Archimedean vector lattice via positive operators.


\begin{lemma}\label{finite nonnegative rank decomposition}
Let $L_1, L_2$ and $L_3$ be bands in a normed lattice $X$ such that $X=L_1\oplus L_2\oplus L_3$. Let $T\colon X\to X$ be a positive operator of the form 
        $$T=\left(\begin{matrix}
             0 & 0 & T_{13}\\
             0 & 0 & 0 \\
             0 & 0 & 0
        \end{matrix}\right)$$ 
        with respect to the decomposition $X=L_1\oplus L_2\oplus L_3$. 
       If $T_{13}$ has a finite nonnegative rank and $\rank^+ (T_{13}) \leq \dim (L_2)$, then there exist positive operators $M, N\colon X\to X$ such that $T=M N$ and 
        $M^2=N^2=N M=0$. 
\end{lemma}

\begin{proof}
Since the nonnegative rank $k:=\rank^+(T_{13})$ is finite, there exist positive operators $R\colon L_3\to \mathbb R^k$ and $L\colon \mathbb R^k \to L_1$ such that $T_{13}=LR$. Clearly, we may assume that $k\geq 1$. 

We claim that $L_2$ contains $k$ pairwise disjoint non-zero positive vectors. If $\dim L_2<\infty$, then for some $n\in \mathbb N_0$ the vector lattice $L_2$ is lattice isomorphic to $\mathbb R^n$ ordered coordinatewise. Since $k\leq n$, $L_2$ clearly contains $k$ pairwise disjoint non-zero positive vectors. 
If $\dim L_2 = \infty$, we use \cite[Theorem 26.10]{LZ71} which states that any infinite-dimensional Archimedean vector lattice contains sets of pairwise disjoint non-zero positive vectors of any finite cardinality. 

Let us now pick any set $\{x_1,\ldots,x_k\}$ in $L_2$ of pairwise disjoint non-zero positive vectors. By $e_i$ we denote the $i$-th standard basis vector of $\mathbb R^k$. The mapping $\iota\colon e_i\mapsto x_i$ $(1\leq i \leq k)$ can be extended to an injective linear operator (denoted again by $\iota$) $\iota\colon \mathbb R^k\to L_2$. Since vectors $x_1,\ldots,x_k$ are pairwise disjoint, $\iota$ is an injective lattice homomorphism. 

Now we are going to construct a positive operator $P\colon L_2\to \mathbb R^k$ such that $P\iota$ is the identity operator on $\mathbb R^k$ as follows. For each $1\leq i\leq k$, by \cite[Theorem 39.3]{Zaa97} there exists a positive bounded linear functional $\phi_i\colon L_2\to \mathbb R$ such that $\phi_i(x_i)=1$. For each $1\leq i\leq k$ we denote by $I_i$ the principal ideal in $L_2$ generated by $x_i$. Applying \cite[Theorem 1.28]{AB06} for the ideal $I_i$, we obtain a positive linear functional $0\leq \varphi_i\leq \phi_i$ on $L_2$ such that $\varphi_i(x_j)=\delta_{ij}$ for $1\leq i,j\leq k$. Since $\phi_i$ is bounded, by \cite[Theorem 25.8]{Zaa97} the functional $\varphi_i$ is bounded as well. If we define $P\colon L_2\to \mathbb R^k$ as $\sum_{i=1}^k e_i\otimes \varphi_i$, then $P$ is a positive operator such that $P\iota$ is the identity operator on $\mathbb R^k$. 

With respect to the decomposition $X=L_1\oplus L_2\oplus L_3$ we define positive operators 
$$ M = \left(\begin{matrix}
    0 & L  P & 0 \\
    0 & 0 & 0 \\
    0 & 0 & 0
        \end{matrix}\right) \qquad \textrm{and} \qquad N=\left(\begin{matrix}
    0 & 0 & 0\\
    0 & 0 & \iota  R \\
    0 & 0 & 0
        \end{matrix}\right)  $$
that satisfy $M^2=N^2=NM=0$ and 
$$ MN = \left(\begin{matrix}
    0 & 0 & LP\iota R \\
    0 & 0 & 0 \\
    0 & 0 & 0
        \end{matrix}\right) =  \left(\begin{matrix}
    0 & 0 & LR \\
    0 & 0 & 0 \\
    0 & 0 & 0
        \end{matrix}\right) =  \left(\begin{matrix}
    0 & 0 & T_{13} \\
    0 & 0 & 0 \\
    0 & 0 & 0
        \end{matrix}\right) = T , 
$$  
which concludes the proof. 
\end{proof}

The following theorem is an infinite-dimensional version of \Cref{main-finite_dimensional}.

\begin{theorem}
For a positive operator $T \colon \ell^p\to \ell^p$ $(1\leq p<\infty)$ the following statements are equivalent. 
\begin{enumerate}
\item There exist positive operators $M,N\colon \ell^p \to \ell^p$ such that $T=M N$ and 
$M^2=N^2=N M=0$. 
\item There exists a positive operator $U\colon \ell^p\to \ell^p$ such that $T=U^2$ and $U^3=0$.
\item There exists a nontrivial band decomposition $\ell^p=L_1\oplus L_2\oplus L_3$ with respect to which the operator $T$ is of the form 
        $$T=\left(\begin{matrix}
                 0 & 0 & T_{13}\\
                 0 & 0 & 0 \\
                 0 & 0 & 0
        \end{matrix}\right)$$ 
        for some positive operator $T_{13}\colon L_3\to L_1$ with $\rank^+ (T_{13}) \leq \dim (L_2)$.
\end{enumerate}
\end{theorem}

\begin{proof}
We may assume that $T \neq 0$.
While the proof of (i)$\Rightarrow$(ii)  can be proved in the same way as (i)$\Rightarrow$(ii) in \Cref{main-finite_dimensional}, the proof of (ii)$\Rightarrow$(iii) requires some explanation. First, using arguments similar to those in \Cref{main-finite_dimensional}, one can find nontrivial bands $L_1, L_2$ and $L_3$ in $\ell^p$ such that $\ell^p=L_1\oplus L_2\oplus L_3$, and positive operators $U_{12}\colon L_2\to L_1$ and $U_{23}\colon L_3\to L_2$ such $T$ is of the form given in (iii) with $T_{13}=U_{12}U_{23} \colon L_3\to L_1$. If $\dim L_2=\infty$, then by definition of the nonnegative rank of a positive operator we have $\rank^+ (T_{13}) \leq \dim (L_2)$. On the other hand, if $\dim L_2<\infty$, since $T_{13}$ factors through the finite-dimensional vector lattice $L_2$, we have $\rank^+ (T_{13}) \leq \dim (L_2)$.

To prove (iii)$\Rightarrow$(i), note first that we can apply Lemma \ref{finite nonnegative rank decomposition} if  $\rank^+ (T_{13})<\infty$. 
Therefore, we may assume that the nonnegative rank of $T_{13}$ is infinite. In particular, $L_2$ is infinite-dimensional.
Since $L_2$ is a band in $\ell^p$, it is lattice isometric to $\ell^p$. Moreover, if $L_1$ were finite-dimensional, then  the following diagram  
\[
\begin{tikzpicture}[>=Stealth, node distance=2.5cm, on grid]
  \node (X) at (0,0) {$L_3$};
  \node (Y) at (4,0) {$L_1$};
  \node (Z) at (2,-2) {$L_1$};

  \draw[->] (X) -- (Y) node[midway, above] {$T_{13}$};
  \draw[->] (X) -- (Z) node[midway, left, xshift=-4pt, yshift=-2pt] {$T_{13}$};
  \draw[->] (Z) -- (Y) node[midway, right, xshift=4pt, yshift=-2pt] {$I$};
\end{tikzpicture}
\] 
would show that $T_{13}\colon L_3\to L_1$ factors through $L_1$. However, this leads to a contradiction with $\rank^+ (T_{13})=\infty$, and we conclude that $\dim L_1=\infty$. By the same reasoning as above, it follows that $L_1$ is also lattice isometric to $\ell^p$.
With respect to the band decomposition 
$\ell^p=L_1\oplus L_2\oplus L_3\cong \ell^p\oplus \ell^p\oplus L_3$ we define positive operators 
$$ M = \left(\begin{matrix}
    0 & I & 0 \\
    0 & 0 & 0 \\
    0 & 0 & 0
        \end{matrix}\right) \qquad \textrm{and} \qquad N=\left(\begin{matrix}
    0 & 0 & 0\\
    0 & 0 & T_{13} \\
    0 & 0 & 0
        \end{matrix}\right).$$
Now it is easy to prove that $T=M N$ and $M^2=N^2=N M=0$. 
\end{proof}

Finally, we consider operators on the Banach lattice $L^p[0,1]$ that has no atoms.

\begin{theorem}
For a positive operator $T\colon L^p[0,1]\to L^p[0,1]$ $(1\leq p<\infty)$ the following statements are equivalent. 
\begin{enumerate}
\item $T=MN-NM$ for some positive operators $M,N\colon L^p[0,1]\to L^p[0,1]$ with $M^2=N^2=0$.
\item There exist positive operators $M,N\colon L^p[0,1]\to L^p[0,1]$ such that $T=MN$ and $M^2=N^2=NM=0$.
\item There exists a positive operator $U\colon L^p[0,1]\to L^p[0,1]$ such that $T=U^2$ and $U^3=0$.
\item There exists a nontrivial band decomposition $L^p[0,1]=B_1\oplus B_2\oplus B_3$ with respect to which the operator $T$ is of the form 
        $$ T=\left(\begin{matrix}
                 0 & 0 & T_{13}\\
                 0 & 0 & 0 \\
                 0 & 0 & 0
        \end{matrix}\right) $$ 
        for some positive operator $T_{13}\colon B_3\to B_1$.
\end{enumerate}
\end{theorem}

  \begin{proof}
  We may assume that the operator $T$ is non-zero. \\
  (i)$\Rightarrow$(iv) Let us write $T=MN-NM$ for some positive square-zero operators $M,N\colon L^p[0,1]\to L^p[0,1]$. Since $0\leq NM\leq MN$, we have 
    $$0\leq (M+N)^3 \leq M^3+3MN^2+3M^2N+N^3=0,$$
    so that $(M+N)^3=0$.  Note that $(M+N)^2 \neq 0$, since otherwise $0 \leq T = M N - N M \leq M N + N M = (M+N)^2 = 0$, and so $T = 0$.
    
    Let us define the bands $B_1=\mathcal N(M+N)$, $B_2=\mathcal N(M+N)^d \cap \mathcal N((M+N)^2)$ and $B_3=\mathcal N((M+N)^2)^d$.
    Since $(M+N)^3 = 0$, the band $B_1$ is non-zero. 
    Since $(M+N)^2 \neq 0$, the band $B_3$ is non-zero as well.  To prove that the band $B_2$ is also non-zero, take a non-zero nonnegative function $g \in B_3$. Then 
    $f = (M+N)g$ is a non-zero nonnegative function in $B_2$.
    
    With respect to the decomposition $L^p[0,1]=B_1\oplus B_2\oplus B_3$ the operators $M+N$ and $(M+N)^2$ are of the form 
    $$M+N=\left(\begin{matrix}
    0 & \star & \star\\
    0 & 0 & \star \\
    0 & 0 & 0
        \end{matrix}\right) \qquad \textrm{and} \qquad (M+N)^2=\left(\begin{matrix}
    0 & 0 & \star\\
    0 & 0 & 0 \\
    0 & 0 & 0
        \end{matrix}\right).$$
    Since $0\leq T=MN-NM\leq MN+NM=(M+N)^2$, the operator $T$ is also of the form 
    $$T=\left(\begin{matrix}
    0 & 0 & \star\\
    0 & 0 & 0 \\
    0 & 0 & 0
        \end{matrix}\right).$$
        
(iv)$\Rightarrow$(ii) Since $L^p[0,1]$ is atomless and $T$ is non-zero, the bands $B_1$, $B_2$ and $B_3$ are non-zero atomless separable Banach lattices. 
    By \cite[Theorem 7.1]{Bohnenblust}, they are all lattice isometric to $L^p[0,1]$. Hence, without any loss of generality we may assume $B_1=B_2=B_3=L^p[0,1]$. With respect to the decomposition $L^p[0,1]=L^p[0,1]\oplus L^p[0,1]\oplus L^p[0,1]$ we define operators $M$ and $N$ as 
    $$M=\left(\begin{matrix}
    0 & T_{13} & 0\\
    0 & 0 & 0 \\
    0 & 0 & 0
    \end{matrix}\right) \qquad \textrm{and} \qquad N=\left(\begin{matrix}
    0 & 0 & 0\\
    0 & 0 & I \\
    0 & 0 & 0
    \end{matrix}\right).$$ 
A direct calculation shows that we have $M^2=N^2=NM=0$ and $MN=T$.

Since it is clear that (ii) implies (i), the assertions (i), (ii) and (iv) are equivalent.

The implication (ii)$\Rightarrow$(iii) can be proved in the same way as (i)$\Rightarrow$(ii) in \Cref{main-finite_dimensional}.
It remains to prove (iii)$\Rightarrow$(iv). Define the bands $B_1=\mathcal N(U)$, $B_2=\mathcal N(U)^d\cap \mathcal N(U^2)$ and $B_3=\mathcal N(U^2)^d$. 
Since $U^3 = 0$, the band $B_1$ is non-zero. Since $U^2 \neq 0$, the band $B_3$ is also non-zero. 
To prove that the band $B_2$ is non-zero as well, take a non-zero nonnegative function $g \in B_3$. Then 
    $f = U g$ is a non-zero nonnegative function in $B_2$.
With respect to the decomposition $\mathbb R^n = B_1\oplus B_2\oplus B_3$ the positive operators $U$ and $U^2$ are of the form 
    $$U=\left(\begin{matrix}
    0 & U_{12} & U_{13}\\
    0 & 0 & U_{23} \\
    0 & 0 & 0
        \end{matrix}\right) \qquad \textrm{and} \qquad T=U^2=\left(\begin{matrix}
    0 &  0 & U_{12} U_{23}\\
    0 & 0 & 0 \\
    0 & 0 & 0
        \end{matrix}\right) , $$
and so $T_{13} = U_{12} U_{23}$.
\end{proof}

\subsection*{Acknowledgements}
The first author was supported by the Slovenian Research and Innovation Agency program P1-0222.
The second author was supported by the Slovenian Research and Innovation Agency program P1-0222 and grant N1-0217.

\end{document}